\documentclass{article}
\usepackage[english]{babel}
\usepackage[utf8]{inputenc}
\usepackage[OT4]{fontenc}
\usepackage{amsmath}
\usepackage{dsfont}
\usepackage{graphicx}
\usepackage{amsthm}
\usepackage{amssymb}

\title{Automatic parametrization and mesh deformation for CFD optimization\\{\sc a technical report}}
\author{\L{}ukasz \L{}aniewski-Wo\l{}\l{}k\footnote{llaniewski@meil.pw.edu.pl, C-CFD group, Warsaw University of Technology. http://c-cfd.meil.pw.edu.pl/}}

\newcommand{\R}{\mathbb{R}}

\begin{document}
\maketitle
{\it Notice: It is a stub of an technical report. It will be extended in future revisions}
\begin{abstract}
We present an automatic and memory efficient methods of morphing-based parametrization of shapes in CFD optimization. Method is based on Kriging and Radial Basis Function interpolation methods.
\end{abstract}
\section{Introduction}
In recent years CFD optimization gained wide recognition as a design technique. With the increase use of adjoint methodologies, shape optimization for fluid flow problems is destined to be a part of industrial-scale product development in coming years. For a survey of CFD optimization methods we refer the reader to the book by Mohammadi and Pironneau\cite{mohammadi_applied_2001,jameson_aerodynamic_2003}.

In all optimization problems it is vital to construct a parametrization of the design space. In many cases the parametrization of shape is achieved by application of morphing (so called morphing box approach) to the surface geometry. In this paper we present a morphing technique in which Radial Basis Functions are used for point displacements and morphing points are selected through a maximum variance criterion.

\section{Method details}

\subsection{Interpolation}

Now let us select a covariance function $\kappa$, which is positive-definite. For example:
\[\kappa(d)=e^{-\frac{d^2}{2\theta^2}}\]

Let us also define function $K$ as $K(x,y)=\kappa(\|x-y\|)$ of covariance between points in $R^3$. Also let us extend the definition to sets of $\R^3$ points:
\[K\left((x_1,x_2,\cdots,x_n),(y_1,y_2,\cdots,y_m)\right)=\left[\begin{array}{cccc}
K(x_1,y_1)& K(x_1,y_2)&&K(x_1,y_m)\\
K(x_2,y_1)& K(x_2,y_2)&&K(x_2,y_m)\\
&&\cdots&\\
K(x_n,y_1)& K(x_n,y_2)&&K(x_n,y_m)\\\end{array}\right]\]

Now let us select a set of $m$ points $M$ which we will call morphing nodes. For a fixed $M$ and a vector of displacements $d$ (matrix of dimension $m\times 3$ we can define a displacement function $\hat{m}:\R^3\rightarrow\R^3$:
\[\hat{m}(x) = d^T K(M,M)^{-1} K(M,x)\]
As the displacement function is defined everywhere in $\R^3$ we can use it for morphing both the surface geometry and volume mesh. It is important to notice that the smoothness of the displacement is dependent solemnly on the smoothness of $\kappa$. For increased control of smoothness the Mattern function can be chosen.

Additionally, basing on Kriging statistical approach, let us define the posteriori variance function $\hat{\sigma^2}:\R^3\rightarrow\R$:
\[\hat{\sigma^2}(x) = \kappa(0) - K(x,M) K(M,M)^{-1} K(M,x)\]
The variance function can be interpreted as a measure of the possible displacement that is ``not parametrized''. It is low in regions where the morphing nodes have high influence and $1$ in regions where there is virtually no displacement done.

\subsection{Algoritm}

Let $P$ be the set of all point in a mesh and $S\subset P$ we the points of the surface mesh. Now we can execute the algorithm of construction of parametrization:

The algorithm of selection of the points is as follows:
\begin{enumerate}
\item Construct the function $\hat{\sigma^2}$ for $M$.\label{step:loop}
\item Find the point $x\in S$ with the maximum $\hat{\sigma^2}$
\item If the maximal value of the variance is below a desired level (or the number of morphing nodes is satisfactory), go to~\ref{step:out}
\item Add $x$ to $M$ and go to~\ref{step:loop}
\item Construct and save matrix $W = K(M,M)^{-1} K(M,P)$\label{step:out}
\end{enumerate}

After this pre-calculation of the parametrization we can calculate the morphing displacement of all the mesh nodes by multiplying a displacement matrix $d$ by the matrix $W$:
\[\hat{m}(P)=d^TW\]

\subsection{Fixing}
In many real world applications, we need to fix parts of the geometry, which cannot be morphed. Presented method can be easily extended to accommodate such need. For a subset $F\subset\R^3$ of the space, we can define a function $d_F(x)$ which is a distance of a point $x$ from the set $F$ (and zero for points inside the set). Now let us define a function $f$ as:
\[f(x) = \kappa(0)-\kappa(d_F(x))\]
We can construct a new $K$ which will take into account that all points in $F$ are fixed:
\[K(x,y)=\kappa(\|x-y\|)f(x)f(y)\]
Resulting $K$ is still positive-define. It is important to notice that this construction works for any positive $f$ which is $0$ on the set of points that are to be fixed, and tends to $1$ otherwise. The function $f$ is constructed basing on $\kappa$ to preserve the influence radius $\theta$ and smoothness of the $\kappa$ across whole parametrization

\section{Applications}
{\it This section will be extended in future revisions}

\section{Summary}
Method for automatic parametrization of shape in CFD optimization cases was showed in detail. The presented approach can be used for both morphing of the surface mesh as the volume mesh, which makes it very useful for application to real-world optimization problems.

\section{Acknowledgements}
This work was supported from 7 Framework Program EU project FLOWHEAD.

\bibliographystyle{elsarticle-num}
\bibliography{../../main}{}

\end{document}